\documentclass[11pt]{amsart}

\usepackage{amsmath,amssymb,amsfonts,amstext,amsthm}
\usepackage{mathtools,latexsym,amscd}
\usepackage{graphicx}
\usepackage[T1]{fontenc}
\usepackage{microtype}
\usepackage{color}
\usepackage{parskip}
\usepackage[pagebackref,colorlinks]{hyperref}
\usepackage{enumitem}
\usepackage{comment}
\usepackage{tikz,tikz-cd}
\usetikzlibrary{decorations.markings}

\newtheorem{thm}{Theorem}[section]
\newtheorem{mainthm}{Theorem}
\newtheorem{lemma}[thm]{Lemma}

\theoremstyle{definition}
\newtheorem{defn}[thm]{Definition}
\newtheorem{prop}[thm]{Proposition}
\newtheorem{coro}[thm]{Corollary}

\newtheorem{re}[thm]{Remark}

\newcommand{\op}{\operatorname}
\newcommand{\ra}{\rightarrow}
\newcommand{\hra}{\hookrightarrow}

\newcommand{\sk}{\hskip}
\newcommand{\mb}{\mathbb}

\newcommand{\ep}{\epsilon}
\newcommand{\al}{\alpha}
\newcommand{\be}{\beta}
\newcommand{\de}{\delta}

\newcommand{\om}{\omega}
\newcommand{\ga}{\gamma}
\newcommand{\Ga}{\Gamma}

\newcommand{\La}{\Lambda}

\newcommand*\nbar[1]{\overline{#1}}

\hypersetup{
  pdftitle={Bounded cohomology of measure-preserving homeomorphism groups of non-orientable surfaces},
  pdfauthor={Michael Brandenbursky and Lior Menashe},
  pdfsubject={Bounded cohomology of transformation groups on non-orientable surfaces},
  pdfkeywords={bounded cohomology, non-orientable surfaces, hyperbolically embedded subgroups}
}

\begin{document}
\raggedbottom

\title[Bounded cohomology of homeomorphism groups]{Bounded cohomology of measure-preserving homeomorphism groups of non-orientable surfaces}
\author{Michael Brandenbursky}\thanks{M.B. was partially supported by the Israel Science Foundation grant 823/23.}
\address{Department of Mathematics, Ben Gurion University, Israel}
\email{brandens@bgu.ac.il} 
\author{Lior Menashe}\thanks{L.M. was partially supported by the Israel Science Foundation grant 823/23.}
\address{Department of Mathematics, Ben Gurion University, Israel}
\email{lior.me10@gmail.com} 

\begin{abstract}
Let $N_g$ be a closed non-orientable surface of genus $g\geq3$, and let
$\op{Homeo}_0(N_g,\mu)$ be the identity component of its group of
measure-preserving homeomorphisms. For every $n\geq2$ we construct a linear
injection
\[
\nbar{EH}_b^n(F_2)\hookrightarrow
\nbar{EH}_b^n\bigl(\op{Homeo}_0(N_g,\mu)\bigr),
\]
where $\nbar{EH}_b^n$ denotes reduced exact bounded cohomology. In
particular, $H_b^2(\op{Homeo}_0(N_g,\mu))$ and
$H_b^3(\op{Homeo}_0(N_g,\mu))$ are infinite-dimensional. The main new
issue is low genus. A natural figure-eight subgroup
$F_2\leq\pi_1(N_g)$ is hyperbolically embedded for $g\geq4$, whereas in
genus three hyperbolic embeddedness fails for every $\pi_1$-injective
figure-eight with orientable regular neighborhood. We overcome this
obstruction using the amalgam decomposition
$\pi_1(N_g)=F_2*_{\mathbb Z}F_{g-2}$ and an isometric bounded-cohomology
extension theorem for graphs of groups with amenable edge groups. Combined
with the Gambaudo--Ghys transfer and annular pushing, this yields the
injection above.
\end{abstract}

\maketitle

\section{Introduction}

Gromov's seminal work \cite{Gromov} established bounded cohomology as a
fundamental tool in geometry and topology, with applications to the geometry
of manifolds, stable commutator length, amenable groups, and symplectic
geometry. In low degrees it is understood for many classes of finitely
generated groups. In higher degrees, however, even the bounded cohomology of
the non-abelian free group $F_2$ remains much less understood; see
\cite{Frigerio}.

Since the early 2000s, bounded cohomology of transformation groups has
received considerable attention. Degree-two bounded classes are closely
related to quasimorphisms, which detect geometric and dynamical phenomena
such as entropy and symplectic spectral invariants; see
\cite{Entropy,G-invariant}. Higher-degree bounded classes of transformation
groups were studied more recently by Brandenbursky--Marcinkowski
\cite{Bounded}, Kimura \cite{Kimura}, and Nitsche \cite{Nitsche}. In the
orientable case, the Gambaudo--Ghys--Polterovich transfer method provides a
powerful mechanism for producing such classes. The non-orientable case has
a genuinely different low-genus feature: the natural geometric free
subgroup needed by the transfer argument is hyperbolically embedded for
$g\geq4$, but in genus three it necessarily fails even the almost
malnormality condition. Resolving this obstruction is one of the main
points of the present paper.

Let $N_g$ be a closed non-orientable surface of genus $g\geq3$. Fix a
hyperbolic metric on $N_g$ and let $\mu$ be its Riemannian area measure,
normalized to have total area one. Equivalently,
the lift of $\mu$ to the orientation double cover
$\Sigma_{g-1}\to N_g$ is induced by the lifted hyperbolic area form. Let
$\op{Homeo}(N_g,\mu)$ be the group of $\mu$-preserving homeomorphisms,
equipped with the compact-open topology, and let
$\op{Homeo}_0(N_g,\mu)$ denote its identity component. 
Fathi's local-contractibility theorem \cite{Fathi} implies that this
identity component is path connected. In particular, every element of
$\op{Homeo}_0(N_g,\mu)$ is joined to the identity by an isotopy through
$\mu$-preserving homeomorphisms. Our main result is the following.

\begin{mainthm}\label{T:main}
Let $g\geq3$. Then
\[
\dim H_b^2(\op{Homeo}_0(N_g,\mu))
=
\dim H_b^3(\op{Homeo}_0(N_g,\mu))
=
\infty.
\]
\end{mainthm}

For every $g\geq3$ and every $n\geq2$ we
construct a linear injection
\begin{equation}\label{eq:intro-injection}
\nbar{EH}_b^n(F_2)\hookrightarrow
\nbar{EH}_b^n\bigl(\op{Homeo}_0(N_g,\mu)\bigr).
\end{equation}
Since $\nbar{EH}_b^2(F_2)$ and $\nbar{EH}_b^3(F_2)$ have dimension
$2^{\aleph_0}$, the theorem follows; in fact, in these two degrees the
targets contain subspaces of dimension $2^{\aleph_0}$.

We now describe the main ingredients. Choose a figure-eight
$i:S^1\vee S^1\hookrightarrow N_g$ contained in an orientable once-holed
torus. Writing $G_g=\pi_1(N_g)$ and identifying the image of $i_*$ with
$H=F_2$, van Kampen gives an amalgam decomposition $G_g=F_2*_{\mathbb Z}F_{g-2}$.
For $g\geq4$ we prove directly that this geometric subgroup $H$ is
hyperbolically embedded in $G_g$ with respect to a natural finite relative
generating set. The proof includes an explicit computation of the relative
metric. We also give an independent proof using Bowditch's relative
hyperbolicity criterion. The extension theorem of
Frigerio--Pozzetti--Sisto \cite{Sisto} then gives surjectivity of restriction
from the reduced exact bounded cohomology of $G_g$ to that of $F_2$.

Genus three is qualitatively different. In the presentation
\[
G_3=\langle a,b,c\mid[a,b]c^2=1\rangle,
\qquad H=\langle a,b\rangle,
\]
one has $\langle c^2\rangle\subseteq H\cap cHc^{-1}$, so $H$ is not
almost malnormal and hence cannot be hyperbolically embedded. We prove a
geometric version of this obstruction: no $\pi_1$-injective embedded
figure-eight in $N_3$ with orientable regular neighborhood has almost
malnormal image. Rather than replacing the figure-eight by a less geometric
subgroup, we use the cyclic amalgam itself. The edge group is amenable, and
the isometric extension theorem of Bucher--Burger--Frigerio--Iozzi--Pagliantini--Pozzetti
\cite{BBFIPP} produces, for every $g\geq3$ and $n\geq2$, an isometric linear
right inverse
\[
J_g^n:H_b^n(F_2)\longrightarrow H_b^n(G_g)
\]
to restriction. This works uniformly in all genera and is stronger than the
surjectivity required for the proof.

Finally, we combine this cohomological section with the
Gambaudo--Ghys--Polterovich transfer. The figure-eight lies in an orientable
once-holed torus, so the annular pushing construction can be carried out locally inside the
non-orientable surface. We prove a norm approximation
showing that, after pullback by the pushing homomorphisms, the transferred
class converges to a positive multiple of its restriction to $F_2$.
Composing the transfer with the isometric section therefore gives the
injection \eqref{eq:intro-injection}. This separates the two genuinely
different issues in the proof: the geometric behavior of the figure-eight
subgroup and the bounded-cohomological extension across the cyclic amalgam.

\begin{re}
Kumar in his thesis \cite{Rishi} constructed an infinite-dimensional space
of homogeneous quasimorphisms on $\op{Homeo}_0(N_g,\mu)$.
Our argument gives a direct bounded-cohomological proof of the
infinite-dimensionality in degree two as well as degree three, and the
construction \eqref{eq:intro-injection} applies in every degree $n\geq2$.
\end{re}

\textbf{Acknowledgments.}
M.B. and L.M. acknowledge the support of Israel Science Foundation grant
823/23. This paper is part of the M.Sc. thesis of the second author. Both
authors thank the Department of Mathematics at Ben Gurion University for
excellent working conditions. We used an LLM for computational exploration, algebraic checks,
and editorial assistance. The mathematical ideas and research directions were developed by us.

\section{Preliminaries}

All cohomology groups in this paper have trivial real coefficients, which we
omit from the notation. Transformation groups are regarded as discrete groups
when taking their bounded cohomology. For $g\geq3$ we write
\[
G_g=\pi_1(N_g,z),\qquad
\mathcal H_g=\op{Homeo}_0(N_g,\mu).
\]
For a homomorphism $f$, we denote by $f^*$ its pullback in bounded cohomology,
and by $f^n$ the induced map in degree $n$ on reduced exact bounded cohomology.
These pullbacks are seminorm non-increasing on bounded cohomology and exact
bounded cohomology, and norm non-increasing on their reduced versions.
Indeed, the cochain pullback is given by precomposition and does not increase
the supremum norm; the inequality passes to cohomology and to the quotient
by zero-seminorm classes.

\subsection{Bounded cohomology}\label{sec:Bounded}

Let $G$ be a group and let $C_b^n(G)$ be the space of bounded functions
$G^{n+1}\to\mb R$. We use the left homogeneous convention: the action is
$(g\cdot f)(g_0,\ldots,g_n)=f(g^{-1}g_0,\ldots,g^{-1}g_n)$,
and $C_b^n(G)^G$ denotes the invariant subspace. The coboundary is
\[
\de^n f(g_0,\ldots,g_{n+1})
=\sum_{j=0}^{n+1}(-1)^j f(g_0,\ldots,\widehat{g_j},\ldots,g_{n+1}).
\]
The cohomology of this homogeneous cochain complex is $H_b^n(G)$.
Equivalently, one can use the inhomogeneous cochain complex
$\bar C_b^n(G)=\ell^\infty(G^n)$, with $G^0$ a singleton. The isometric
identification with homogeneous cochains is given by
\[
\phi(f)(g_1,\ldots,g_n)
=f(1,g_1,g_1g_2,\ldots,g_1\cdots g_n),
\]
and the inhomogeneous differential is $\bar\de=\phi\de\phi^{-1}$.
We refer to \cite{Frigerio} for these conventions and the basic properties
of bounded cohomology.

\subsubsection{Exact and reduced bounded cohomology}
The exact bounded cohomology is
\[
EH_b^n(G)=\ker\bigl(c_G^n:H_b^n(G)\longrightarrow H^n(G)\bigr),
\]
where $c_G^n$ is the comparison map. The supremum norm on cochains induces
\[
\|\xi\|=\inf\{\|c\|_\infty:[c]=\xi\}
\]
on $H_b^n(G)$. Write $N^n(G)$ for its subspace of zero-seminorm classes.
We use the notation
\[
\nbar{H}_b^n(G)=H_b^n(G)/N^n(G),\qquad
\nbar{EH}_b^n(G)=EH_b^n(G)/(EH_b^n(G)\cap N^n(G)).
\]
Both reduced spaces carry genuine norms. Naturality of the comparison map
and the norm inequality for pullbacks imply that every group homomorphism
induces maps on these spaces.

\subsection{Quasimorphisms}
A map $f:G\to\mb R$ is a \textit{quasimorphism} if
\[
D(f):=\sup_{g,h\in G}|f(gh)-f(g)-f(h)|<\infty.
\]
It is \textit{homogeneous} if $f(g^k)=kf(g)$ for all $g\in G$ and
$k\in\mb Z$. Let $Q(G,\mb R)$ and $Q^h(G,\mb R)$ denote the corresponding
vector spaces. The map $f\mapsto[\bar\de f]$ gives
\[
EH_b^2(G)\cong
\frac{Q(G,\mb R)}{\bar C_b^1(G)+\op{Hom}(G,\mb R)}
\cong\frac{Q^h(G,\mb R)}{\op{Hom}(G,\mb R)};
\]
see \cite[Chapter 2]{Frigerio}.

\subsection{Hyperbolically embedded subgroups}
We recall the definition from \cite[Section 4]{Hyperbolically}.
Cayley graphs have edges of length one, and paths may traverse generator
edges in either direction. Let $H\leq G$ and $X\subseteq G$ with
$G=\langle X\cup H\rangle$. The alphabet $X\sqcup H$ is a disjoint union,
so labels from its two parts are distinguished even when they represent
the same group element. A path in $\Ga(G,X\sqcup H)$ is
\textit{admissible} if it uses no edge of the distinguished subgraph
$\Ga(H,H)$. In particular, $H$-labeled edges in other left cosets of $H$
are allowed. For $h,k\in H$, define $\hat d(h,k)$ as the minimum length of
an admissible path from $h$ to $k$, with value $\infty$ when no such path
exists. This is an extended metric on $H$.

\begin{defn}\label{d:hyp_emb_subgroup}
We write $H\hra_h(G,X)$ if:
\begin{enumerate}[label=(\alph*)]
\item $G=\langle X\cup H\rangle$;
\item $\Ga(G,X\sqcup H)$ is hyperbolic;
\item every ball of finite radius in $(H,\hat d)$ is finite.
\end{enumerate}
We write $H\hra_hG$ if such an $X$ exists.
\end{defn}

We will verify these conditions directly. The following edge-addition
criterion will be used for condition (b).
\begin{prop}[{\cite[Lemma 5.5]{Acylindrically}; see also
\cite[Corollary 2.4]{Hyperbolicity}}]\label{p:edge-addition}
Suppose that a graph $\Gamma_H$ is obtained from a hyperbolic graph $\Ga$ by
adding edges. If there exists $M>0$ such that every $\Ga$-geodesic joining
the endpoints of an edge of $\Gamma_H$ has diameter at most $M$ in $\Gamma_H$, then
$\Gamma_H$ is hyperbolic.
\end{prop}

\begin{defn}
A subgroup $H\leq G$ is \textit{malnormal} if
$H\cap xHx^{-1}=\{1\}$ for every $x\in G\setminus H$. It is
\textit{almost malnormal} if all these intersections are finite.
\end{defn}

\begin{prop}[{\cite[Proposition 4.33]{Hyperbolically}}]\label{p:almost-malnormal}
A hyperbolically embedded subgroup is almost malnormal.
\end{prop}

\subsubsection{Examples}
\begin{itemize}
\item The subgroup $\langle a^2,ab,ab^{-1}\rangle\leq F(a,b)$ is not
hyperbolically embedded. It is the index-two kernel of the map sending
both $a$ and $b$ to the nontrivial element of $\mb Z/2\mb Z$, and hence is
not almost malnormal.
\item The subgroup $\langle a^nba^{-n}:n\geq0\rangle\leq F(a,b)$ has
infinite intersection with its conjugate by $a$, although $a$ does not
belong to it. It is therefore not hyperbolically embedded.
\item For $k>2$, the free factor $F(a_1,a_2)$ is hyperbolically embedded in
$F(a_1,\ldots,a_k)$ relative to $X=\{a_3,\ldots,a_k\}$. The ordinary
geodesic between vertices of a left coset of this free factor stays in
that coset, so Proposition~\ref{p:edge-addition} applies. Normal forms for the
free product show that $\hat d(1,h)=\infty$ for $h\neq1$, which proves
condition (c).
\end{itemize}

For the geometric figure-eight embedding used below, hyperbolic embeddedness
holds for $g\geq4$, but fails for $g=3$. In genus three we instead use an
isometric bounded-cohomology extension over an amenable amalgam; see
Proposition~\ref{p:isometric-section}.

\section{Gambaudo--Ghys map}

We recall the transfer construction of \cite{Bounded}, extending the
constructions of Gambaudo--Ghys \cite{Gambaudo} and Polterovich
\cite{Polterovich}, with conventions compatible with our left homogeneous
cochains. Fix $g\geq3$. Outside the cut locus of $z$, let $\al_{z,x}$ be the
unique minimizing geodesic from $z$ to $x$; on the cut locus, which has
measure zero, choose the paths arbitrarily. This gives a measurable family
of reference paths, as in \cite[Section 2]{Bounded}. Put
$\al_{x,z}=\al_{z,x}^{-1}$. For $f\in\mathcal H_g$, choose an isotopy
$\{f_t\}_{0\leq t\leq1}$ from the identity to $f$ and set
\[
\ga(f,x)=[\al_{z,x}*\{f_t(x)\}_{t\in[0,1]}*\al_{f(x),z}]\in G_g.
\]
Here $*$ denotes traversal from left to right, while multiplication of
based loop classes is taken in composition order: $[q][p]$ means first
$p$, then $q$. With this convention the trajectory cocycle has the
identity stated in Proposition~\ref{p:ga_fh} below.

\begin{prop}
The element $\ga(f,x)$ is independent of the chosen isotopy.
\end{prop}
\begin{proof}
For a loop $\{u_t\}$ in $\mathcal H_g$ and a based loop $\de$ in $N_g$,
the square $(s,t)\mapsto u_t(\de(s))$ shows that the evaluation loop
$t\mapsto u_t(z)$ commutes with $[\de]$. Thus the image of
\[
(ev_z)_*:\pi_1(\mathcal H_g,\op{Id})\longrightarrow G_g
\]
is contained in $Z(G_g)$. This center is trivial: a central deck
transformation of $\mb H^2$ would commute with two hyperbolic elements
with distinct pairs of fixed points at infinity, and consequently fix at
least three boundary points. It must therefore be the identity.
Such elements exist because $G_g$ is a cocompact, non-elementary surface
group; see \cite{Bridson,Harpe}. The same argument applies at any
basepoint $x$. Concatenating one isotopy with the reverse of another
therefore gives a null-homotopic evaluation loop, proving the claim.
\end{proof}

\begin{prop}\label{p:ga_fh}
For $f,h\in\mathcal H_g$ and $x\in N_g$,
\[
\ga(fh,x)=\ga(f,h(x))\ga(h,x).
\]
\end{prop}
\begin{proof}
Use the isotopy that first follows $h_t$ and then $f_t\circ h$.
The trajectory of $x$ is the trajectory of $h_t(x)$ followed by that
of $f_t(h(x))$. The two intermediate reference paths cancel, and the
formula follows from the multiplication convention.
\end{proof}

To use left homogeneous cochains, define
\begin{equation}\label{eq:eta}
\eta(f,x)=\ga(f^{-1},x)^{-1}.
\end{equation}
Proposition~\ref{p:ga_fh} gives
\begin{equation}\label{eq:eta-equivariance}
\eta(hf,x)=\eta(h,x)\eta(f,h^{-1}(x)).
\end{equation}
For $c\in C_b^n(G_g)^{G_g}$ put
\begin{equation}\label{eq:transfer-cochain}
I_b^n(c)(f_0,\ldots,f_n)
=\int_{N_g}c\bigl(\eta(f_0,x),\ldots,\eta(f_n,x)\bigr)\,d\mu(x).
\end{equation}
The trajectory cocycle is measurable; see \cite[Section 2]{Bounded} and
\cite[Section 3.B]{Volume}. Thus the integral is well defined.

\begin{prop}\label{p:I_b}
The maps $I_b^n$ take homogeneous cochains to homogeneous cochains,
commute with the coboundary, and satisfy
\[
\|I_b^n(c)\|_\infty\leq\mu(N_g)\|c\|_\infty.
\]
\end{prop}
\begin{proof}
The norm estimate is immediate, and the alternating sum defining the
coboundary commutes with integration. For homogeneity, use
\eqref{eq:eta-equivariance}, left invariance of $c$, and invariance of
$\mu$ to obtain
\begin{align*}
I_b^n(c)(hf_0,\ldots,hf_n)
&=\int_{N_g}c\bigl(\eta(f_0,h^{-1}(x)),\ldots,
                         \eta(f_n,h^{-1}(x))\bigr)\,d\mu(x)\\
&=I_b^n(c)(f_0,\ldots,f_n).\qedhere
\end{align*}
\end{proof}

It follows that $I_b^n$ induces
$\Ga_b^n:H_b^n(G_g)\to H_b^n(\mathcal H_g)$. It also preserves exact
classes. Indeed, for a fixed $f$, the map $x\mapsto\ga(f,x)$ has
essentially finite image, as in \cite[Lemma 2.2]{Bounded}.
Consequently \eqref{eq:transfer-cochain} is also defined on arbitrary
homogeneous cochains. The same calculation as in Proposition~\ref{p:I_b}
therefore gives a cochain map on ordinary homogeneous cochains, and the
resulting square with the comparison maps commutes. The norm estimate then
gives the induced map
\[
\nbar{E\Ga}_b^n:
\nbar{EH}_b^n(G_g)\longrightarrow
\nbar{EH}_b^n(\mathcal H_g).
\]

\section{Proofs}

The reduced exact bounded cohomology of $F_2$ has dimension
$2^{\aleph_0}$ in degrees two and three; for degree two see
\cite{Brooks} and \cite[Section 2.6]{Frigerio}, and for degree three see
\cite{Soma} and \cite[Corollary 3]{Sisto}. We will find an embedded figure-eight
$i:S^1\vee S^1\hookrightarrow N_g$, contained in an orientable once-holed torus,
for which
\[
i^n:\nbar{EH}_b^n(G_g)\longrightarrow\nbar{EH}_b^n(F_2)
\]
is surjective for every $n\geq2$. For $g\geq4$ we first prove this using
a direct hyperbolic-embedding argument. For $g=3$, where that geometric
subgroup is not hyperbolically embedded, we construct an isometric
linear right inverse to restriction using an amenable amalgam.
The latter construction in fact works for every $g\geq3$.

As in \cite{Bounded}, annular pushing maps then give homomorphisms
$\rho_\ep:F_2\to\mathcal H_g$ and a constant $\La>0$ such that
\[
\bigl\|\rho_\ep^n\nbar{E\Ga}_b^n(\xi)-\La i^n(\xi)\bigr\|
\longrightarrow0
\quad\text{as }\ep\to0
\]
for every $\xi\in\nbar{EH}_b^n(G_g)$. Thus the diagram
\[
\begin{tikzcd}[column sep=large,row sep=large]
\nbar{EH}_b^n(G_g)
  \arrow[r,"\nbar{E\Ga}_b^n"]\arrow[d,"i^n"']
&\nbar{EH}_b^n(\mathcal H_g)\arrow[dl,"\rho_\ep^n"]\\
\nbar{EH}_b^n(F_2)&
\end{tikzcd}
\]
commutes up to multiplication by $\La$ and an error tending to zero.
This will prove
\[
\dim H_b^n(\mathcal H_g)\geq
\dim\nbar{EH}_b^n(F_2),\qquad n\geq2.
\]

\subsection{An elementary retraction when \texorpdfstring{$g\geq5$}{g >= 5}}\label{sec:genus-ge5}
We retain the elementary argument in higher genus, but use $j$ and $r$
to distinguish this figure-eight from the figure-eight embedding $i$ used later.

\begin{prop}\label{p:retraction}
For $g\geq5$ there are an embedding $j:S^1\vee S^1\hookrightarrow N_g$
and a retraction $r:N_g\to S^1\vee S^1$ such that $r\circ j=\op{Id}$.
\end{prop}
\begin{proof}
Write $N_g=\Sigma_2\#N_{g-4}$. Choose an embedded figure-eight given by
one meridian in each of the two handles, joined at the common basepoint
as in Figure~\ref{fig:retraction}. In the standard CW description the
surface group has presentation
\[
\left\langle x_1,y_1,x_2,y_2,d_1,\ldots,d_{g-4}\ \middle|\
[x_1,y_1][x_2,y_2]d_1^2\cdots d_{g-4}^2=1\right\rangle.
\]
Define a map of the one-skeleton to the chosen figure-eight that is the
identity on $x_1,x_2$ and collapses the remaining generator loops.
The attaching loop of the two-cell maps to a null-homotopic loop, so
by the usual CW extension criterion this map extends over the two-cell. It
is the required retraction. Equivalently, its induced homomorphism sends
$x_1,x_2$ to the two free generators and all other generators to $1$.
\end{proof}

\begin{figure}[htbp]
    \centering

    \tikzset{every picture/.style={line width=0.75pt}}

    \begin{tikzpicture}[x=0.75pt,y=0.75pt,yscale=-1,xscale=1]

        \draw  [line width=1.5]  (336.54,207.77) -- (341.27,215.67) -- (346.96,208.43) ;
        \draw [line width=1.5]    (130.03,237.82) .. controls (120.79,225.11) and (122.1,205.13) .. (130.3,192.25) ;
        \draw [line width=1.5]    (130.3,192.25) .. controls (145.71,165.48) and (197.86,156.28) .. (242.48,180.97) ;
        \draw [line width=1.5]    (130.03,237.82) .. controls (142.27,259.44) and (203.11,270.38) .. (243.56,241) ;
        \draw [line width=1.5]    (242.48,180.97) .. controls (273.04,158.04) and (330.05,163.32) .. (350.46,183.63) ;
        \draw [line width=1.5]    (243.56,241) .. controls (251.47,247.12) and (309.11,271.01) .. (349.56,241.62) ;
        \draw [line width=1.5]    (350.46,183.63) .. controls (334.03,211.38) and (343.61,238.53) .. (350.58,240.64) ;
        \draw [line width=1.5]    (350.46,183.63) .. controls (360.33,191.78) and (364.84,223.86) .. (350.58,240.64) ;
        \draw  [color={rgb, 255:red, 0; green, 0; blue, 0 },draw opacity=1][line width=1.5]  (349.66,183.62) .. controls (349.66,183.07) and (350.12,182.63) .. (350.67,182.63) .. controls (351.22,182.64) and (351.66,183.1) .. (351.66,183.65) .. controls (351.65,184.2) and (351.19,184.64) .. (350.64,184.63) .. controls (350.09,184.63) and (349.65,184.17) .. (349.66,183.62) -- cycle ;
        \draw  [color={rgb, 255:red, 0; green, 0; blue, 0 },draw opacity=1][line width=1.5]  (349.66,183.62) -- (351.66,183.65) ;
        \draw  [color={rgb, 255:red, 0; green, 0; blue, 0 },draw opacity=1][line width=1.5]  (350.67,182.63) -- (350.64,184.63) ;
        \draw  [color={rgb, 255:red, 0; green, 0; blue, 0 },draw opacity=1][line width=1.5]  (348.78,240.61) .. controls (348.79,240.06) and (349.24,239.62) .. (349.8,239.63) .. controls (350.35,239.64) and (350.79,240.09) .. (350.78,240.64) .. controls (350.77,241.19) and (350.32,241.64) .. (349.76,241.63) .. controls (349.21,241.62) and (348.77,241.16) .. (348.78,240.61) -- cycle ;
        \draw  [color={rgb, 255:red, 0; green, 0; blue, 0 },draw opacity=1][line width=1.5]  (348.78,240.61) -- (350.78,240.64) ;
        \draw  [color={rgb, 255:red, 0; green, 0; blue, 0 },draw opacity=1][line width=1.5]  (349.8,239.63) -- (349.76,241.63) ;
        \draw  [line width=1.5]  (354.07,208.69) -- (359.74,215.95) -- (364.49,208.06) ;
        \draw  [color={rgb, 255:red, 0; green, 0; blue, 0 },draw opacity=1][line width=0.75]  (241.68,180.96) .. controls (241.69,180.41) and (242.15,179.97) .. (242.7,179.98) .. controls (243.25,179.98) and (243.69,180.44) .. (243.68,180.99) .. controls (243.68,181.54) and (243.22,181.98) .. (242.67,181.98) .. controls (242.12,181.97) and (241.68,181.51) .. (241.68,180.96) -- cycle ;
        \draw  [color={rgb, 255:red, 0; green, 0; blue, 0 },draw opacity=1][line width=0.75]  (241.68,180.96) -- (243.68,180.99) ;
        \draw  [color={rgb, 255:red, 0; green, 0; blue, 0 },draw opacity=1][line width=0.75]  (242.7,179.98) -- (242.67,181.98) ;
        \draw  [color={rgb, 255:red, 0; green, 0; blue, 0 },draw opacity=1][line width=0.75]  (242.76,240.98) .. controls (242.77,240.43) and (243.23,239.99) .. (243.78,240) .. controls (244.33,240.01) and (244.77,240.46) .. (244.76,241.01) .. controls (244.75,241.57) and (244.3,242.01) .. (243.75,242) .. controls (243.19,241.99) and (242.75,241.54) .. (242.76,240.98) -- cycle ;
        \draw  [color={rgb, 255:red, 0; green, 0; blue, 0 },draw opacity=1][line width=0.75]  (242.76,240.98) -- (244.76,241.01) ;
        \draw  [color={rgb, 255:red, 0; green, 0; blue, 0 },draw opacity=1][line width=0.75]  (243.78,240) -- (243.75,242) ;
        \draw [line width=1.5]    (161,211.72) .. controls (174.82,223.94) and (201.86,221.35) .. (216.04,209.57) ;
        \draw [line width=1.5]    (172.12,217.3) .. controls (177.12,203.97) and (201.15,202.34) .. (206.95,215.43) ;
        \draw [line width=1.5]    (267.02,211.35) .. controls (280.83,223.57) and (304.86,221.94) .. (319.04,210.15) ;
        \draw [line width=1.5]    (277.94,216.52) .. controls (282.13,204.58) and (304.14,203.92) .. (309.97,215.01) ;
        \draw [color={rgb, 255:red, 40; green, 70; blue, 140 },draw opacity=1][line width=0.75]    (150.33,203.56) .. controls (166.9,174.81) and (218.47,181.6) .. (242.22,210.97) ;
        \draw [color={rgb, 255:red, 40; green, 70; blue, 140 },draw opacity=1][line width=0.75]    (242.22,210.97) .. controls (271.88,177.42) and (314.73,178.68) .. (328.46,196.29) ;
        \draw [color={rgb, 255:red, 40; green, 70; blue, 140 },draw opacity=1][line width=0.75]    (153.93,229.62) .. controls (150.38,226.16) and (144.11,217.87) .. (150.33,203.56) ;
        \draw [color={rgb, 255:red, 40; green, 70; blue, 140 },draw opacity=1]   (153.93,229.62) .. controls (170.89,247.17) and (217.8,238.2) .. (242.22,210.97) ;
        \draw [color={rgb, 255:red, 40; green, 70; blue, 140 },draw opacity=1]   (242.22,210.97) .. controls (267.56,241.36) and (318.96,243.44) .. (331.14,217.34) ;
        \draw [color={rgb, 255:red, 40; green, 70; blue, 140 },draw opacity=1]   (328.46,196.29) .. controls (333.37,202.37) and (334.24,210.39) .. (331.14,217.34) ;

    \end{tikzpicture}

    \caption{The embedding $j:S^1\vee S^1\hookrightarrow N_5$ and the retraction $r:N_5\to S^1\vee S^1$.}  
    \label{fig:retraction}

\end{figure}

\begin{coro}
For $n\geq2$, $j^n\circ r^n=\op{Id}$. In particular, $j^n$ is surjective
and $r^n$ is an isometric injection.
\end{coro}
\begin{proof}
Functoriality gives the identity. Since both pullbacks are norm
non-increasing, $\|u\|=\|j^nr^n(u)\|\leq\|r^n(u)\|\leq\|u\|$.
\end{proof}

\subsection{The free subgroup and its amalgam decomposition}\label{sec:free-subgroup}
For the rest of the proof, let $g\geq3$ and put $m=g-2$. Decompose
\[
N_g=T\cup_\de N_{m,1},\qquad
T\cong\Sigma_{1,1},\quad
\de=\partial T=\partial N_{m,1}.
\]
Here $\Sigma_{1,1}$ is a compact once-holed torus, and $N_{m,1}$ is a
compact non-orientable surface of genus $m$ with one boundary component.
Choose meridian and longitude loops $\al,\be$ meeting transversely only
at $z$, whose union is a deformation retract of $T$. This defines
$i:S^1\vee S^1\hookrightarrow N_g$ in the orientable part of
$\mb T^2\#N_m$; see Figure~\ref{fig:figure-eight}. Write $a,b$ for their
based homotopy classes and $w=[a,b]=aba^{-1}b^{-1}$.

With suitable choices of paths to the boundary and of generator
orientations, van Kampen's theorem gives
\begin{equation}\label{eq:surface-amalgam}
\begin{split}
G_g&=\left\langle a,b,c_1,\ldots,c_m\ \middle|\
[a,b]c_1^2\cdots c_m^2=1\right\rangle\\
&=H*_{C}B,\qquad H=F(a,b),\quad B=F(c_1,\ldots,c_m),
\end{split}
\end{equation}
where, putting $v=c_1^2\cdots c_m^2$, the abstract edge group and its
injections are
\begin{equation}\label{eq:edge-inclusions}
C_0=\langle t\rangle\cong\mb Z,\qquad
H\xleftarrow{\ t\mapsto w\ }C_0
 \xrightarrow{\ t\mapsto v^{-1}\ }B.
\end{equation}
Both boundary words are nontrivial and of infinite order in their free
groups, so the maps in \eqref{eq:edge-inclusions} are injective.
The normal-form theorem for amalgams \cite[Chapter I, Section 1]{Serre}
therefore identifies $H$ and $B$ with subgroups of $G_g$, and
\[
H\cap B=C=\langle w\rangle=\langle v\rangle,\qquad w=v^{-1}.
\]
Thus $i_*:F_2\to G_g$ is injective for every $g\geq3$. We identify $F_2$
with $H$ by this map and retain $i^n$ for the induced restriction on
reduced exact bounded cohomology.

\begin{figure}[htbp]
    \centering

    \tikzset{every picture/.style={line width=0.75pt}}

    \begin{tikzpicture}[x=0.75pt,y=0.75pt,yscale=-1,xscale=1]

    \draw  [line width=1.5]  (151.28,90.11) -- (156.52,97.68) -- (161.72,90.08) ;
    \draw [line width=1.5]    (165,63) .. controls (150,85.5) and (157,115) .. (165.97,121.18) ;
    \draw [line width=1.5]    (165,63) .. controls (176,68) and (180.23,104.39) .. (165.97,121.18) ;
    \draw  [color={rgb, 255:red, 0; green, 0; blue, 0 },draw opacity=1][line width=1.5]  (164.66,63.62) .. controls (164.66,63.07) and (165.12,62.63) .. (165.67,62.63) .. controls (166.22,62.64) and (166.66,63.1) .. (166.66,63.65) .. controls (166.65,64.2) and (166.19,64.64) .. (165.64,64.63) .. controls (165.09,64.63) and (164.65,64.17) .. (164.66,63.62) -- cycle ;
    \draw  [color={rgb, 255:red, 0; green, 0; blue, 0 },draw opacity=1][line width=1.5]  (164.66,63.62) -- (166.66,63.65) ;
    \draw  [color={rgb, 255:red, 0; green, 0; blue, 0 },draw opacity=1][line width=1.5]  (165.67,62.63) -- (165.64,64.63) ;
    \draw  [line width=1.5]  (169.34,90.31) -- (174.46,97.95) -- (179.77,90.43) ;
    \draw  [color={rgb, 255:red, 40; green, 70; blue, 140 },draw opacity=1][line width=0.75]  (96.68,118.96) .. controls (96.69,118.41) and (97.15,117.97) .. (97.7,117.98) .. controls (98.25,117.98) and (98.69,118.44) .. (98.68,118.99) .. controls (98.68,119.54) and (98.22,119.98) .. (97.67,119.98) .. controls (97.12,119.97) and (96.68,119.51) .. (96.68,118.96) -- cycle ;
    \draw  [color={rgb, 255:red, 40; green, 70; blue, 140 },draw opacity=1][line width=0.75]  (96.68,118.96) -- (98.68,118.99) ;
    \draw  [color={rgb, 255:red, 40; green, 70; blue, 140 },draw opacity=1][line width=0.75]  (97.7,117.98) -- (97.67,119.98) ;
    \draw [line width=1.5]    (46,120) .. controls (37,109) and (34,86.5) .. (44.3,68.25) ;
    \draw [line width=1.5]    (44.3,68.25) .. controls (65,34.5) and (132,35) .. (165,63) ;
    \draw [line width=1.5]    (86.12,96.3) .. controls (91.12,82.97) and (115.15,81.34) .. (120.95,94.43) ;
    \draw  [color={rgb, 255:red, 0; green, 0; blue, 0 },draw opacity=1][line width=1.5]  (164.78,120.61) .. controls (164.79,120.06) and (165.24,119.62) .. (165.8,119.63) .. controls (166.35,119.64) and (166.79,120.09) .. (166.78,120.64) .. controls (166.77,121.19) and (166.32,121.64) .. (165.76,121.63) .. controls (165.21,121.62) and (164.77,121.16) .. (164.78,120.61) -- cycle ;
    \draw  [color={rgb, 255:red, 0; green, 0; blue, 0 },draw opacity=1][line width=1.5]  (164.78,120.61) -- (166.78,120.64) ;
    \draw  [color={rgb, 255:red, 0; green, 0; blue, 0 },draw opacity=1][line width=1.5]  (165.8,119.63) -- (165.76,121.63) ;
    \draw  [color={rgb, 255:red, 40; green, 70; blue, 140 },draw opacity=1] (60,91.71) .. controls (60,76.4) and (79.03,64) .. (102.5,64) .. controls (125.97,64) and (145,76.4) .. (145,91.71) .. controls (145,107.01) and (125.97,119.41) .. (102.5,119.41) .. controls (79.03,119.41) and (60,107.01) .. (60,91.71) -- cycle ;
    \draw [dashed] [color={rgb, 255:red, 40; green, 70; blue, 140 },draw opacity=1]   (106,140.67) .. controls (121,149) and (127,92.33) .. (105,98.33) ;
    \draw [color={rgb, 255:red, 40; green, 70; blue, 140 },draw opacity=1]   (106,140.67) .. controls (96.2,137) and (94.2,108) .. (105,98.33) ;
    \draw [line width=1.5]    (46,120) .. controls (63,145) and (129.2,151) .. (165.97,121.18) ;
    \draw [line width=1.5]    (76,90.72) .. controls (89.82,102.94) and (116.86,100.35) .. (131.04,88.57) ;
    \end{tikzpicture}
    \caption{The figure-eight embedding $i:S^1\vee S^1\hookrightarrow N_g$, illustrated for $g=3$. The complementary piece has $g-2$ crosscaps.}
    \label{fig:figure-eight}        
\end{figure}

\subsection{A direct hyperbolic embedding for \texorpdfstring{$g\geq4$}{g >= 4}}\label{sec:direct}

\begin{prop}\label{p:F2_hyp_emb}
Let $g\geq4$, $m=g-2$, and use the notation of
\eqref{eq:surface-amalgam}. Then
\[
H\hra_h(G_g,X),\qquad X=\{c_1,\ldots,c_m\}.
\]
More precisely, for the relative metric defined by this finite set $X$,
\begin{equation}\label{eq:relative-metric}
\hat d(1,h)=
\begin{cases}
2m|k|,&h=w^k,\quad k\in\mb Z,\\
\infty,&h\in H\setminus C.
\end{cases}
\end{equation}
\end{prop}

\begin{proof}
We verify Definition~\ref{d:hyp_emb_subgroup} directly. Condition (a) is
immediate. Write
\[
\Ga=\Ga(G_g,\{a,b,c_1,\ldots,c_m\}),\qquad
\Gamma_H=\Ga(G_g,X\sqcup H).
\]
The graph $\Ga$ is hyperbolic because $G_g$ acts properly and cocompactly
on $\mb H^2$; see \cite[Chapters I.8 and III.H]{Bridson}.
We first establish the path estimates needed for (b) and (c).

\textit{Subpaths leaving and returning to a factor.}
Let $\mathcal T$ be the Bass--Serre tree of $H*_CB$. Its vertices are
$G_g/H\sqcup G_g/B$ and its edges are $G_g/C$, with $uC$ joining
$uH$ to $uB$; see \cite[Chapter I, Section 4]{Serre}.
Map a group vertex $u$ to the midpoint of $uC$. Map an $H$-labelled edge
through the corresponding $H$-vertex, and a $c_j^{\pm1}$-edge through
the corresponding $B$-vertex. The same construction applies to
$\Ga$, with its $a,b$-edges mapped through $H$-vertices.

Consequently, a path that leaves $B$ and returns to $B$, with endpoints
$x,y\in B$ and all intermediate vertices outside $B$, has $xC=yC$: its image avoids the vertex $B$ and
must return along the same incident edge of the tree. The analogous
statement holds for $H$.

\textit{Overlap of boundary axes.}
Let $T_B=\Ga(B,X)$ be the ordinary Cayley tree of $B$, and let $L$ be
the axis of $v=c_1^2\cdots c_m^2$, containing $1$. Its positive labels
repeat the word $c_1c_1c_2c_2\cdots c_mc_m$.
Each oriented two-letter subword in this cyclic word occurs in exactly
one position modulo its length $2m$: these subwords are $c_jc_j$ and
$c_jc_{j+1}$, where $c_{m+1}=c_1$. Here we use $m\geq2$.
If $L$ and $uL$ share two consecutive edges, their positive orientations
agree on the common segment, since all positive labels belong to $X$
and none belongs to $X^{-1}$. The two-letter word then determines the
same phase on both lines. Uniqueness of labelled paths in a Cayley tree
forces the entire lines, with their phases, to agree. Their phase-zero
vertices are respectively $C$ and $uC$, so $u\in C$. Therefore
\begin{equation}\label{eq:axis-overlap}
\op{diam}_{T_B}(L\cap uL)\leq1\qquad(u\in B\setminus C).
\end{equation}
If $p_L:T_B\to L$ is nearest-point projection, it follows that
\begin{equation}\label{eq:coset-projection}
\op{diam}_{T_B}p_L(uC)\leq1\qquad(u\in B\setminus C).
\end{equation}
Indeed, $uC\subset uL$; two disjoint lines in a tree project to a single
point, and intersecting lines project onto their intersection. See
Figure~\ref{fig:direct-embedding}(b) for the case of a one-edge intersection.

\textit{A weighted normal-form observation.}
Give the $a,b$-edges of $\Ga$ length $1$ and the $c_j$-edges length
$2/m$, and write $d_*$ for the resulting path metric. The intrinsic
weighted lengths of the two representatives of every edge-group element
agree:
\begin{equation}\label{eq:weighted-boundary}
|w^k|_{\{a,b\}}=4|k|
=\frac{2}{m}|v^{-k}|_X.
\end{equation}
Both factors are isometrically embedded on vertices for these weighted
metrics. We give the normal-form argument explicitly. Let $W$ be a word in
$\{a^{\pm1},b^{\pm1},c_1^{\pm1},\ldots,c_m^{\pm1}\}$ representing an
element $h\in H$. Split $W$ into maximal alternating $H$- and $B$-blocks
and replace each block by a shortest word in its factor. If more than one
block remains, the amalgam normal-form theorem implies that some block
representing an element of $C$ can be transferred to the adjacent factor;
otherwise the alternating word would be reduced and could not represent
an element of $H$. Transfer such a block using
\eqref{eq:weighted-boundary} and merge the adjacent blocks. This does not
increase weighted length and decreases the number of blocks. Iterating
leaves a single block. If that block lies in $B$, then, because it represents
the element $h\in H$, it belongs to $H\cap B=C$ and can be transferred once
more to $H$ without changing weighted length. We have therefore produced an
$H$-word of weighted length at most that of $W$. The reverse inequality is
immediate because an intrinsic $H$-word is also an ambient word. The argument
for an element of $B$ is identical. Thus both factors are isometrically
embedded. In particular,
\begin{equation}\label{eq:weighted-isometry}
d_*(1,h)=|h|_{\{a,b\}}\ (h\in H),\qquad
d_*(x,y)=\frac{2}{m}d_{T_B}(x,y)\ (x,y\in B).
\end{equation}

\textit{Condition (b): hyperbolicity of $\Gamma_H$.}
Suppose first that $m>2$. If an ordinary $\Ga$-geodesic $q$ with
endpoints $1,h\in H$ used a $c_j$-edge, its weighted length would be
strictly smaller than its ordinary length. Since an $H$-word also joins
these endpoints, \eqref{eq:weighted-isometry} would imply
\[
|h|_{\{a,b\}}=d_*(1,h)
\leq\ell_*(q)<\ell(q)=d_\Ga(1,h)\leq|h|_{\{a,b\}},
\]
a contradiction. Hence every such geodesic stays in the ordinary
$H$-Cayley tree. Its vertices have diameter at most $1$ in $\Gamma_H$;
including points in edge interiors gives a bound of $2$.

It remains to consider $m=2$, i.e. $g=4$. Now all the weights are one,
so \eqref{eq:weighted-isometry} says that the ordinary factor Cayley
trees are isometrically embedded on vertices in $\Ga$.
Let $q$ be a subpath of a $\Ga$-geodesic whose endpoints lie in $H$
and whose intermediate vertices lie outside $H$. After left translation
by its initial vertex, the preceding Bass--Serre-tree observation gives
endpoints $1,v^k$, and no intermediate vertex of $q$ lies in $H$.

We claim that $q$ stays in $B$. Otherwise, consider each maximal
subpath $r$ of $q$ whose endpoints $x,y$ lie in $B$ and whose intermediate
vertices lie outside $B$, and replace $r$ by the unique $T_B$-geodesic
joining $x$ to $y$. The preceding Bass--Serre-tree observation gives
$xC=yC$,
so $x,y\in uC$ for some $u\in B$. These endpoints are intermediate
vertices of $q$, so $u\notin C$; such a subpath cannot begin or end at
an endpoint of $q$, because then its first or last edge would be an
$a$- or $b$-edge in $H$, which is impossible. Since $r$ is a subpath
of a geodesic and $B$ is isometrically embedded on vertices, each
replacement preserves length.
The resulting path is a $T_B$-geodesic from $1$ to $v^k$, hence is the
segment in $L$. Each nontrivial replacement segment lies in both $L$
and $uL$ and has length at least $4$, because its distinct endpoints
belong to $uC$. This contradicts \eqref{eq:axis-overlap}.
Thus $q$ is the segment in the adjacent $B$-tree between its endpoints.

Successive vertices of $C$ on $L$ are at distance $4$, so every point
of this segment is within distance $2$ of $H$. It follows, by decomposing
any $\Ga$-geodesic with endpoints in $H$ into its maximal subpaths whose
endpoints lie in $H$ and whose intermediate vertices lie outside $H$,
together with the parts lying in $H$, that the whole geodesic is within
distance $2$ of the vertex set $H$. Since $H$ has diameter one in $\Gamma_H$, its
diameter in $\Gamma_H$ is at most $2+1+2=5$.
Figure~\ref{fig:direct-embedding}(a) illustrates such a genus-four
subpath and the added edge in the base coset.

In both cases, left translation gives a uniform bound for every
$\Ga$-geodesic joining the endpoints of an added $H$-edge. For the
original generator edges the required bound is immediate.
Proposition~\ref{p:edge-addition} proves that $\Gamma_H$ is hyperbolic.

\begin{figure}[htbp]
\centering
\definecolor{proofblue}{RGB}{40,70,140}
\begin{tikzpicture}[x=1cm,y=1cm,font=\footnotesize,
  line width=0.75pt,
  directed/.style={postaction={decorate},
    decoration={markings,mark=at position .55 with {\arrow{>}}}}]
\node[anchor=west,font=\small] at (0,2.0)
  {(a) $g=4$: the two boundary words have length $4$.};
\coordinate (a0) at (.7,0);
\coordinate (a1) at (3,1.0);
\coordinate (a2) at (5.2,1.3);
\coordinate (a3) at (7.4,1.0);
\coordinate (a4) at (9.7,0);
\coordinate (b1) at (3,-.85);
\coordinate (b2) at (5.2,-1.1);
\coordinate (b3) at (7.4,-.85);
\draw[proofblue,thick,directed] (a0)--node[above,sloped] {$c_2^{-1}$} (a1);
\draw[proofblue,thick,directed] (a1)--node[above] {$c_2^{-1}$} (a2);
\draw[proofblue,thick,directed] (a2)--node[above] {$c_1^{-1}$} (a3);
\draw[proofblue,thick,directed] (a3)--node[above,sloped] {$c_1^{-1}$} (a4);
\draw[directed] (a0)--node[below,sloped] {$a$} (b1);
\draw[directed] (b1)--node[below] {$b$} (b2);
\draw[directed] (b2)--node[below] {$a^{-1}$} (b3);
\draw[directed] (b3)--node[below,sloped] {$b^{-1}$} (a4);
\draw[densely dashed] (a0)--node[above] {$H$-edge in $\Gamma_H$}
  node[below] {length $1$} (a4);
\foreach \p in {a0,a4,b1,b2,b3}
  \fill (\p) circle (1.5pt);
\foreach \p in {a1,a2,a3}
  \draw[proofblue,fill=white] (\p) circle (1.6pt);
\node[left] at (a0) {$h_0$};
\node[right] at (a4) {$h_0w$};
\node[above,proofblue] at (5.2,1.42) {$q\subset h_0L$};
\node[below] at (5.2,-1.48) {the lower path lies in $H$};
\end{tikzpicture}

\smallskip
\begin{tikzpicture}[x=1cm,y=1cm,font=\footnotesize,line width=0.75pt]
\node[anchor=west,font=\small] at (0,2.4)
  {(b) $g\geq4$: an allowed $H$-edge cannot give a long projection.};
\coordinate (x) at (1.7,1.2);
\coordinate (y) at (8.7,1.2);
\coordinate (px) at (4.7,0);
\coordinate (py) at (5.7,0);
\draw[<->] (.25,0)--(10.5,0);
\node[right] at (10.5,0) {$L$};
\draw[proofblue,thick,<->] (.4,1.6)--(x)--(3.5,.7)--(px)--(py)
  --(6.8,.7)--(y)--(10,1.6);
\node[right,proofblue] at (10,1.6) {$uL$};
\draw[line width=1.6pt] (px)--(py);
\draw[proofblue,densely dashed]
  (x).. controls (3.6,2.05) and (6.8,2.05)..(y);
\node[above,proofblue] at (5.2,1.78) {allowed $H$-edge in $uH$};
\fill[proofblue] (x) circle (1.7pt);
\fill[proofblue] (y) circle (1.7pt);
\fill (px) circle (1.5pt);
\fill (py) circle (1.5pt);
\node[below left] at (x) {$x\in uC$};
\node[below right] at (y) {$y\in uC$};
\node[below left] at (px) {$p_L(x)$};
\node[below right] at (py) {$p_L(y)$};
\draw[|<->|] (4.7,-.58)--node[below] {$\leq1$} (5.7,-.58);
\end{tikzpicture}
\caption{The estimates in Proposition~\ref{p:F2_hyp_emb}.
(a) Here $h_0\in H$ and $w=[a,b]=c_2^{-2}c_1^{-2}$.
Every point of the upper path is within distance $2$ of $H$,
whose diameter in $\Gamma_H$ is one. The dashed edge lies in the
forbidden base coset $H$.
(b) For $u\in B\setminus C$, an allowed $H$-edge with endpoints in
$uC$ projects to a segment of $L$ of length at most one. The axes are
drawn schematically, not to scale.}
\label{fig:direct-embedding}
\end{figure}

\textit{Condition (c): the relative metric.}
Let $q$ be an admissible path from $1$ to $h\in H$. Its image in
$\mathcal T$ avoids the vertex $H$, since precisely the forbidden
$H$-edges in the base coset would pass through that vertex. The
midpoints of $C$ and $hC$ must therefore lie in the same component of
$\mathcal T\setminus\{H\}$. Hence $hC=C$, or $h\in C$. Thus
\[
\hat d(1,h)=\infty\quad(h\in H\setminus C).
\]

Take a shortest admissible path from $1$ to $v^k$. Such a path exists,
since the word $v^k$ in $X^{\pm1}$ is admissible. It cannot contain a
maximal subpath whose endpoints $x,y$ lie in $B$ and whose intermediate
vertices lie outside $B$. Indeed, the preceding Bass--Serre-tree observation
gives $xC=yC$. The first edge of such a subpath has an $H$-label, so
admissibility implies $x\notin H$. If $x\neq y$, the single $H$-edge
joining $x$ to $y$ lies in $xH\neq H$ and is allowed, giving a shorter
path. If $x=y$, delete the subpath. This is again a contradiction.
Thus a shortest admissible path stays in $B$.

It may still use $H$-edges in other cosets; these must not be discarded.
Every such edge has endpoints in $uC$ for some $u\in B\setminus C$,
because $B\cap H=C$. By \eqref{eq:coset-projection} their projections
to $L$ are at distance at most one; see
Figure~\ref{fig:direct-embedding}(b). The same bound holds for the
endpoints of each $X^{\pm1}$-edge, since projection onto a subtree is
$1$-Lipschitz. Projecting the path edge by edge therefore gives
\[
2m|k|=d_{T_B}(1,v^k)\leq\ell(q).
\]
The $X$-word $v^k$ realizes equality. Since $w=v^{-1}$, this proves
\eqref{eq:relative-metric}. In particular, for every $r\geq0$ the closed
relative ball is
\begin{equation}\label{eq:relative-balls}
\overline B_{\hat d}(1,r)
=\{w^k:k\in\mb Z,\ |k|\leq\lfloor r/(2m)\rfloor\}.
\end{equation}
Left translation by $H$ preserves admissibility, so every finite-radius
ball is finite. This proves condition (c) and completes the proof.
\end{proof}

\begin{re}[Alternative proof using Bowditch's criterion]\label{r:Bowditch}
Here is a second proof of Proposition~\ref{p:F2_hyp_emb}, without the
relative-metric computation. We will show that $H$ is a proper infinite
quasiconvex malnormal subgroup of the non-elementary hyperbolic group
$G_g$. Bowditch's theorem \cite[Theorem 7.11]{Bowditch}, applied to the
conjugacy class of $H$, then shows that $G_g$ is hyperbolic relative to
$H$. By \cite[Proposition 4.28]{Hyperbolically}, this is equivalent to
hyperbolic embeddedness with respect to any finite relative generating
set, in particular $X=\{c_1,\ldots,c_m\}$.

We check the hypotheses. For $m\geq2$, the cyclic word
$v=c_1c_1\cdots c_mc_m$ has no proper period: the cyclic two-letter word
$c_1c_2$ occurs exactly once. Hence $v$ is not a proper power, and $C$
is a maximal cyclic subgroup of the free group $B$, therefore malnormal
in $B$. For completeness, if $C\cap uCu^{-1}$ contains a nontrivial
element, the axes of powers of $v$ and of $uvu^{-1}$ coincide. The
stabilizer of this axis in a free group is generated by the primitive
root of $v$, which here is $v$ itself. Thus $u\in C$; this is also a
standard consequence of the free action on the Cayley tree
\cite[Chapter I]{Serre}.

The subgroup $H$ is then malnormal in $H*_CB$. Indeed, for $x\notin H$,
an element of $H\cap xHx^{-1}$ fixes the segment joining the distinct
vertices $H$ and $xH$ of the Bass--Serre tree. This segment contains a
$B$-vertex with two distinct incident edges. Their stabilizers are, up
to conjugacy, $C$ and $uCu^{-1}$ with $u\in B\setminus C$, and their
intersection is trivial.

Finally, $H$ is quasiconvex. For $g\geq4$, the curve $\de=\partial T$
bounds neither a disk nor a M\"obius band, so it can be represented by a
simple closed geodesic in a hyperbolic metric on $N_g$. After an isotopy,
$T$ therefore has geodesic boundary. A component $\widetilde T$ of its
preimage in $\mb H^2$ is an intersection of hyperbolic half-planes, hence
convex. Its stabilizer is $H$, and $\widetilde T/H=T$ is compact. Thus
an $H$-orbit is coarsely dense in $\widetilde T$ and quasiconvex in
$\mb H^2$. The Milnor--\v Svarc lemma and stability of quasi-geodesics
transfer this to a Cayley graph of $G_g$; see
\cite[Chapters I.8 and III.H]{Bridson}. The subgroup is infinite and
proper by the amalgam normal form. Bowditch's criterion now applies.
\end{re}

The consequence for bounded cohomology follows from the following
extension result.
\begin{thm}[{Frigerio--Pozzetti--Sisto, \cite[Proposition 6.1]{Sisto}}]
\label{t:Sisto}
If $H\hra_hG$, restriction $EH_b^n(G)\to EH_b^n(H)$ is surjective for
every $n\geq2$. Consequently, restriction is also surjective on reduced
exact bounded cohomology.
\end{thm}
Indeed, restriction is norm non-increasing, so it sends zero-seminorm
classes to zero-seminorm classes. Hence the surjection on exact bounded
cohomology descends to a surjection on the reduced exact quotients.

\begin{coro}\label{c:restriction-ge4}
For the figure-eight embedding $i$ and $g\geq4$, the map
$i^n:\nbar{EH}_b^n(G_g)\to\nbar{EH}_b^n(F_2)$ is surjective
for every $n\geq2$.
\end{coro}

\subsection{Genus three and isometric extension in bounded cohomology}
\label{sec:g3}

When $g=3$, we have $m=1$, so the complementary surface $N_{m,1}$ is a M\"obius band and
\begin{equation}\label{eq:g3-amalgam}
G_3=\langle a,b,c\mid[a,b]c^2=1\rangle
=F(a,b)*_{\langle[a,b]\rangle=\langle c^{-2}\rangle}\langle c\rangle.
\end{equation}
The edge group remains infinite cyclic. Its inclusion into the
M\"obius-band vertex group is $\mb Z\to\mb Z$, $k\mapsto-2k$;
its image has index two, not order two.

\begin{prop}\label{p:g3-obstruction}
The free subgroup $H=\langle a,b\rangle$ in \eqref{eq:g3-amalgam}
is not hyperbolically embedded in $G_3$ with respect to any relative
generating set.
\end{prop}
\begin{proof}
The orientation character $G_3\to\{\pm1\}$ is trivial on $H$ and
sends $c$ to $-1$, so $c\notin H$. But $c^2=[a,b]^{-1}\in H$, and
\[
\langle c^2\rangle\subseteq H\cap cHc^{-1}.
\]
This intersection is infinite, so Proposition~\ref{p:almost-malnormal}
applies. Directly, in $\Ga(G_3,\{c\}\sqcup H)$ the path
\[
1\xrightarrow{\ c\ }c
 \xrightarrow{\ c^{2k}\in H\ }c^{2k+1}
 \xrightarrow{\ c^{-1}\ }c^{2k}
\]
is admissible for every $k\geq1$: its middle edge lies in $cH$. 
Hence $\hat d(1,c^{2k})\leq3$ and the
relative metric is not proper. This is exactly the shortcut that the
axis-overlap estimate prevents when $g\geq4$.
\end{proof}

\begin{re}\label{r:g3-geometric-obstruction}
In fact, no $\pi_1$-injective embedded figure-eight in $N_3$ with
orientable regular neighborhood has almost malnormal image.
To see this, let $R$ be that regular neighborhood. Since $\chi(R)=-1$,
it is either a once-holed torus or a pair of pants. The complement has
Euler characteristic zero. It has no disk component, since a disk
would kill a nontrivial boundary element of $\pi_1(R)$. Each component
has boundary and nonpositive Euler characteristic; consequently each
is an annulus or a M\"obius band. There are either one or three boundary
components in total, so at least one complementary component is a
M\"obius band. Its core, based using a path in $R$, gives an
orientation-reversing element $x$ outside $\pi_1(R)$ with
$x^2\in\pi_1(R)\setminus\{1\}$. Thus
$\langle x^2\rangle\subseteq\pi_1(R)\cap x\pi_1(R)x^{-1}$ is infinite.
The obstruction is therefore geometric, not merely a poor choice of
the generators $a,b$.
\end{re}

The obstruction does not prevent extension of bounded cohomology. We
use the following specialization of the graph-of-groups theorem of
Bucher, Burger, Frigerio, Iozzi, Pagliantini and Pozzetti.

\begin{thm}[{\cite[Theorem 1.1]{BBFIPP}}]\label{t:amenable-amalgam}
Let $A,B$ be countable groups, let $D$ be amenable, and suppose the
maps $D\to A$ and $D\to B$ are injective. Put $G=A*_DB$, and let
$\iota_A:A\hookrightarrow G$ and $\iota_B:B\hookrightarrow G$ be the
canonical inclusions. For $n\geq2$, write
\[
\begin{aligned}
\operatorname{res}_A=\iota_A^*&:H_b^n(G)\longrightarrow H_b^n(A),\\
\operatorname{res}_B=\iota_B^*&:H_b^n(G)\longrightarrow H_b^n(B)
\end{aligned}
\]
for the induced restriction maps. There is a linear map
\[
\Theta^n:H_b^n(A)\oplus H_b^n(B)\longrightarrow H_b^n(G)
\]
such that, for every $u\in H_b^n(A)$ and $v\in H_b^n(B)$,
\[
\begin{gathered}
\operatorname{res}_A\bigl(\Theta^n(u,v)\bigr)=u,\qquad
\operatorname{res}_B\bigl(\Theta^n(u,v)\bigr)=v,\\
\|\Theta^n(u,v)\|=\max\{\|u\|,\|v\|\}.
\end{gathered}
\]
Thus $\Theta^n$ is an isometric linear embedding for the maximum
seminorm on the direct sum. The first two identities say that it is a
right inverse to the simultaneous restriction map
\[
H_b^n(G)\longrightarrow H_b^n(A)\oplus H_b^n(B),\qquad
\xi\longmapsto\bigl(\operatorname{res}_A(\xi),
                       \operatorname{res}_B(\xi)\bigr).
\]
\end{thm}

No malnormality or maximality assumption on the images of $D$ occurs
in this theorem. In particular, the index-two boundary inclusion of
the M\"obius band is allowed.

\begin{prop}\label{p:isometric-section}
For every $g\geq3$ and $n\geq2$, restriction for the figure-eight embedding
admits an isometric linear right inverse
\[
J_g^n:H_b^n(F_2)\longrightarrow H_b^n(G_g),\qquad
 i^*J_g^n=\operatorname{Id}.
\]
It descends to an isometric linear map
\begin{equation}\label{eq:reduced-section}
\overline J_g^n:\nbar{EH}_b^n(F_2)
 \longrightarrow\nbar{EH}_b^n(G_g),\qquad
 i^n\overline J_g^n=\operatorname{Id}.
\end{equation}
Thus the required restriction $i^n$ is surjective also for $g=3$.
\end{prop}
\begin{proof}
In \eqref{eq:surface-amalgam}, the two vertex groups $H=F_2$ and
$B=F_m$ are countable. The edge group $C_0\cong\mb Z$ is amenable,
and both edge maps are injective by \eqref{eq:edge-inclusions}.
Apply Theorem~\ref{t:amenable-amalgam} to these edge maps, with
$A=H=F_2$ and $G=G_g$. The restriction $\operatorname{res}_A$ in that
theorem is precisely $i^*$. Define
\[
J_g^n(u)=\Theta^n(u,0).
\]
Then $i^*J_g^n(u)=u$. Moreover,
\[
\|u\|=\|i^*J_g^n(u)\|\leq\|J_g^n(u)\|\leq\|(u,0)\|=\|u\|,
\]
so $J_g^n$ is isometric. When $g=3$, $B=\mb Z$ and
$H_b^n(B)=0$ for $n\geq1$, but the argument needs no modification.

We check exactness before passing to reduced cohomology. The surface
$N_g$ is aspherical for $g\geq3$, so it is a $K(G_g,1)$.
Since $N_g$ is closed and non-orientable,
$H^2(N_g;\mb R)=0$, and its cohomology in degrees above two vanishes
by dimension. Hence
\[
H^n(G_g;\mb R)=0\qquad(n\geq2).
\]
Also $H^n(F_2;\mb R)=0$ for $n\geq2$, since a graph is a $K(F_2,1)$.
These are the standard surface and graph cohomology calculations;
see \cite[Sections 1.B and 3.1]{Hatcher}. Therefore
\[
EH_b^n(G_g)=H_b^n(G_g),\qquad EH_b^n(F_2)=H_b^n(F_2)
\quad(n\geq2).
\]
Isometry shows that $J_g^n$ sends zero-seminorm classes to zero-seminorm
classes, so it descends to the map in \eqref{eq:reduced-section}. If
$\eta$ has zero seminorm, then $\|\xi+\eta\|=\|\xi\|$ by the triangle
inequality in both directions; hence the quotient norm is the seminorm of
any representative. It follows that the descended map is again isometric.
The identity with restriction also descends, proving surjectivity.
\end{proof}

\begin{re}
The maps $J_g^n$ and $\overline J_g^n$ are cohomological sections. They
are not asserted to be pullbacks of a group or surface retraction.
In particular, Proposition~\ref{p:isometric-section} does not assert a
hyperbolic embedding in genus three. It supplies exactly the restriction
surjectivity needed below, and gives an alternative to
Corollary~\ref{c:restriction-ge4} in every genus $g\geq4$ as well.
\end{re}

\subsection{Construction of the annular pushing homomorphisms}
\label{sec:pushing}

Fix the figure-eight embedding $i$ from Section~\ref{sec:free-subgroup}. Its two loops
have annular neighborhoods inside the orientable subsurface $T$.
Thus the construction of \cite[Lemma 3.1]{Bounded} applies locally,
regardless of the orientability of $N_g$. We provide the details,
keeping the cochain conventions of Section~3.

\begin{lemma}\label{L:main}
For every $g\geq3$ there exist $\La>0$ and homomorphisms
$\rho_\ep:F_2\to\mathcal H_g$, $0<\ep<1/2$, such that, for every
$n\geq2$ and $\xi\in\nbar{EH}_b^n(G_g)$,
\begin{equation}\label{eq:approximation}
\bigl\|\rho_\ep^n\nbar{E\Ga}_b^n(\xi)-\La i^n(\xi)\bigr\|
\longrightarrow0\qquad(\ep\to0).
\end{equation}
The constant $\La$ and the family $\rho_\ep$ are independent of $n$
and $\xi$.
\end{lemma}

\begin{proof}
Let $\al,\be$ be the two loops in the image of the figure-eight embedding $i$.
Choose closed smooth
annular neighborhoods $N(\al),N(\be)\subset\operatorname{int}(T)$
whose intersection is a disk containing $z$ in its interior. Put
\[
\La=\mu(N(\al)\cap N(\be))>0.
\]
On either annulus the measure $\mu$ is represented by a smooth positive
area form after choosing an orientation.

Use the model annulus
\[
N=\{(r,\theta):1\leq r\leq2,\ \theta\in\mb R/2\pi\mb Z\},
\qquad\kappa=r\,dr\wedge d\theta.
\]
There is an orientation-preserving diffeomorphism
$n_\al:N(\al)\to N$ for which
\[
n_\al^*(s_\al\kappa)=\nu_\al,\qquad
s_\al=\frac{\mu(N(\al))}{\int_N\kappa}>0,
\]
where $\nu_\al$ represents $\mu$ on $N(\al)$.
Indeed, start with an orientation-preserving parametrization and
apply Moser's argument after this normalization of total area; the
version for manifolds with boundary is due to Banyaga
\cite{Moser,Banyaga}. The parametrization can be chosen so that
increasing $\theta$ represents the chosen orientation of $\al$.
Choose $n_\be$ similarly. Notice that equal areas of the two annuli,
or equality with the Euclidean area of $N$, are not required.

Choose a smooth function $f_\ep:[1,2]\to[0,1]$ that equals $1$ on
$[1+\ep,2-\ep]$ and vanishes in neighborhoods of $1$ and $2$.
The isotopy
\[
P_\ep^t(r,\theta)=(r,\theta+2\pi t f_\ep(r)),\qquad0\leq t\leq1,
\]
preserves $\kappa$ and every constant multiple of it. Conjugate by
$n_\al$, respectively $n_\be$, and extend by the identity outside the
corresponding annulus. This gives measure-preserving isotopies
$P_\ep^t(\al),P_\ep^t(\be)$ of $N_g$. Identifying $F_2$ with $F(a,b)$, define the homomorphism
$\rho_\ep:F_2\to\mathcal H_g$ by
\[
\rho_\ep(a)=P_\ep^1(\al),\qquad
\rho_\ep(b)=P_\ep^1(\be).
\]
Freeness ensures that these assignments define a homomorphism;
injectivity of $\rho_\ep$ is not needed.

Put
\[
A_\ep(\al)=n_\al^{-1}\{1+\ep\leq r\leq2-\ep\},\qquad
B_\ep(\al)=N(\al)\setminus A_\ep(\al),
\]
and define the analogous sets for $\be$. The time-one map
$\rho_\ep(a)$ fixes $A_\ep(\al)$ pointwise, although the isotopy traces
one full loop around $\al$ there. Consequently, for
$x\in A_\ep(\al)$, $\ga(\rho_\ep(a),x)$ is conjugate to $i_*(a)$;
for $x\notin N(\al)$ it is trivial. The analogous statements hold for
$b$; see Figure~\ref{fig:annular-pushing}.

\begin{figure}[htbp]
    \centering

    \tikzset{every picture/.style={line width=0.75pt}}

    \begin{tikzpicture}[x=0.75pt,y=0.75pt,yscale=-1,xscale=1]
        \draw  [line width=1.5]  (75.03,89.87) .. controls (85.71,44.05) and (117.75,12.38) .. (146.59,19.15) .. controls (175.43,25.91) and (190.15,68.54) .. (179.47,114.37) .. controls (168.79,160.19) and (136.75,191.85) .. (107.91,185.09) .. controls (79.07,178.32) and (64.35,135.69) .. (75.03,89.87) -- cycle ;
        \draw  [line width=1.5]  (109.32,99.51) .. controls (113.61,81.11) and (125.2,68.1) .. (135.22,70.45) .. controls (145.23,72.8) and (149.88,89.62) .. (145.59,108.01) .. controls (141.3,126.41) and (129.71,139.41) .. (119.69,137.06) .. controls (109.68,134.72) and (105.03,117.9) .. (109.32,99.51) -- cycle ;
        \draw  [line width=1.5]  (253.17,104.06) .. controls (253.17,58.57) and (289.92,21.71) .. (335.26,21.71) .. controls (380.6,21.71) and (417.36,58.57) .. (417.36,104.06) .. controls (417.36,149.54) and (380.6,186.4) .. (335.26,186.4) .. controls (289.92,186.4) and (253.17,149.54) .. (253.17,104.06) -- cycle ;
        \draw  [line width=1.5]  (310.99,104.06) .. controls (310.99,90.61) and (321.86,79.71) .. (335.26,79.71) .. controls (348.67,79.71) and (359.53,90.61) .. (359.53,104.06) .. controls (359.53,117.5) and (348.67,128.4) .. (335.26,128.4) .. controls (321.86,128.4) and (310.99,117.5) .. (310.99,104.06) -- cycle ;
        \draw   (116.11,154.15) .. controls (98.67,150.29) and (89.51,124.48) .. (95.67,96.5) .. controls (101.83,68.52) and (120.96,48.97) .. (138.41,52.83) .. controls (155.85,56.7) and (165,82.51) .. (158.85,110.49) .. controls (152.69,138.47) and (133.56,158.02) .. (116.11,154.15) -- cycle ;
        \draw    (253.5,103.5) -- (418,104.65) ;
        \draw [line width=1.5]    (370.04,101.3) -- (370.04,107.91) ;
        \draw [line width=1.5]    (407.36,101.3) -- (407.36,107.91) ;
        \draw  [color={rgb, 255:red, 0; green, 0; blue, 0 },draw opacity=1][line width=1.5]  (369.19,56.52) .. controls (369.2,55.91) and (369.7,55.43) .. (370.31,55.44) .. controls (370.91,55.45) and (371.4,55.95) .. (371.39,56.56) .. controls (371.38,57.16) and (370.88,57.65) .. (370.27,57.64) .. controls (369.67,57.63) and (369.18,57.13) .. (369.19,56.52) -- cycle ;
        \draw  [color={rgb, 255:red, 0; green, 0; blue, 0 },draw opacity=1][line width=1.5]  (369.19,56.52) -- (371.39,56.56) ;
        \draw  [color={rgb, 255:red, 0; green, 0; blue, 0 },draw opacity=1][line width=1.5]  (370.31,55.44) -- (370.27,57.64) ;
        \draw   (277.19,103.06) .. controls (277.19,70.88) and (303.19,44.8) .. (335.26,44.8) .. controls (367.34,44.8) and (393.34,70.88) .. (393.34,103.06) .. controls (393.34,135.23) and (367.34,161.31) .. (335.26,161.31) .. controls (303.19,161.31) and (277.19,135.23) .. (277.19,103.06) -- cycle ;
        \draw  [color={rgb, 255:red, 40; green, 70; blue, 140 },draw opacity=1] (143.37,34.15) .. controls (166.41,39.38) and (178.13,74.51) .. (169.54,112.62) .. controls (160.94,150.72) and (135.29,177.37) .. (112.25,172.14) .. controls (89.2,166.9) and (77.48,131.77) .. (86.07,93.67) .. controls (94.67,55.57) and (120.32,28.92) .. (143.37,34.15) -- cycle ;
        \draw [color={rgb, 255:red, 40; green, 70; blue, 140 },draw opacity=1]   (143,52.33) .. controls (145,44.33) and (150,33.33) .. (151,36) ;
        \draw   (360.92,159.92) -- (354.09,158.29) -- (358.24,152.63) ;
        \draw  [color={rgb, 255:red, 0; green, 0; blue, 0 },draw opacity=1][line width=1.5]  (148.29,35.92) .. controls (148.3,35.31) and (148.8,34.82) .. (149.4,34.83) .. controls (150.01,34.84) and (150.49,35.34) .. (150.48,35.95) .. controls (150.47,36.56) and (149.98,37.04) .. (149.37,37.03) .. controls (148.76,37.02) and (148.28,36.52) .. (148.29,35.92) -- cycle ;
        \draw  [color={rgb, 255:red, 0; green, 0; blue, 0 },draw opacity=1][line width=1.5]  (148.29,35.92) -- (150.48,35.95) ;
        \draw  [color={rgb, 255:red, 0; green, 0; blue, 0 },draw opacity=1][line width=1.5]  (149.4,34.83) -- (149.37,37.03) ;
        \draw  [color={rgb, 255:red, 0; green, 0; blue, 0 },draw opacity=1][line width=1.5]  (141.29,53.92) .. controls (141.3,53.31) and (141.8,52.82) .. (142.4,52.83) .. controls (143.01,52.84) and (143.49,53.34) .. (143.48,53.95) .. controls (143.47,54.56) and (142.98,55.04) .. (142.37,55.03) .. controls (141.76,55.02) and (141.28,54.52) .. (141.29,53.92) -- cycle ;
        \draw  [color={rgb, 255:red, 0; green, 0; blue, 0 },draw opacity=1][line width=1.5]  (141.29,53.92) -- (143.48,53.95) ;
        \draw  [color={rgb, 255:red, 0; green, 0; blue, 0 },draw opacity=1][line width=1.5]  (142.4,52.83) -- (142.37,55.03) ;
        \draw (209,100.75) node [anchor=north west][inner sep=0.75pt]   [align=left] {$\ra$};
        \draw (208,82.5) node [anchor=north west][inner sep=0.75pt]  [font=\small] [align=left] {$n_\al$};
        \draw (116.25,142.2) node [anchor=north west][inner sep=0.75pt]  [font=\footnotesize] [align=left] {$\al$};
        \draw (135.5,59) node [anchor=north west][inner sep=0.75pt]  [font=\scriptsize] [align=left] {$z$};
        \draw (359.5,109) node [anchor=north west][inner sep=0.75pt]  [font=\scriptsize] [align=left] {$1$+$\ep$};
        \draw (396,109) node [anchor=north west][inner sep=0.75pt]  [font=\scriptsize] [align=left] {$2\sk 0.125em$-$\sk 0.1em\ep$};
        \draw (342.75,61) node [anchor=north west][inner sep=0.75pt]  [font=\scriptsize] [align=left] {$n_{\al}(x)$};
        \draw (147,41.85) node [anchor=north west][inner sep=0.75pt]  [font=\scriptsize] [align=left] {$x$};
        \draw (303.85,163.5) node [anchor=north west][inner sep=0.75pt]  [font=\scriptsize] [align=left] {$P_{\ep}^{t}(n_{\al}(x))$};
        \draw (148.5,187.5) node [anchor=north west][inner sep=0.75pt]   [align=left] {$N(\al)$};
        \draw (386.96,187.25) node [anchor=north west][inner sep=0.75pt]   [align=left] {$N$};
        \draw (111.75,159.75) node [anchor=north west][inner sep=0.75pt]  [font=\scriptsize] [align=left] {$\ga$};

    \end{tikzpicture}

    \caption{The annular pushing isotopy. For $x\in A_\ep(\al)$, the time-one map fixes $x$ and the trajectory class $\ga(\rho_\ep(a),x)$ is conjugate to $i_*(a)$.}
    \label{fig:annular-pushing}
\end{figure}

Up to boundary sets of measure zero, decompose $N_g$ into
\begin{align*}
A_\ep&=A_\ep(\al)\cap A_\ep(\be),&
A_\ep^a&=A_\ep(\al)\setminus N(\be),\\
A_\ep^b&=A_\ep(\be)\setminus N(\al),&
B_\ep&=B_\ep(\al)\cup B_\ep(\be),\\
D&=N_g\setminus(N(\al)\cup N(\be)).&&
\end{align*}
The five sets displayed above form a measurable partition up to null
boundary sets. Both time-one generators fix $N_g\setminus B_\ep$
pointwise; since they are bijections, $B_\ep$ is invariant as well.
Thus every piece of the partition is $\rho_\ep(F_2)$-invariant. Let
$h_a:F_2\to\langle a\rangle$ send $a$ to $a$ and $b$ to $1$, and define
$h_b$ similarly. Because the intersection
$N(\al)\cap N(\be)$ is a disk, the conjugating paths for the two
trajectory loops can be chosen simultaneously there. Thus, for
$\om\in F_2$, Proposition~\ref{p:ga_fh} gives
\begin{equation}\label{eq:trajectory-pieces}
\ga(\rho_\ep(\om),x)=
\begin{cases}
u_x i_*(\om)u_x^{-1},&x\in A_\ep,\\
u_{x,a} i_*h_a(\om)u_{x,a}^{-1},&x\in A_\ep^a,\\
u_{x,b} i_*h_b(\om)u_{x,b}^{-1},&x\in A_\ep^b,\\
1,&x\in D,
\end{cases}
\end{equation}
where the conjugators depend on $x$, but not on $\om$.
For example, on $A_\ep$ choose, measurably in $x$, a path from $z$ to
$x$ inside the intersection disk and compare it with the reference path
$\al_{z,x}$. On $A_\ep^a$ and $A_\ep^b$ use paths inside the
corresponding annulus. The resulting path classes, and hence the
conjugators, are measurable and take values in the countable group $G_g$.
Since these points are fixed by $\rho_\ep(F_2)$,
\eqref{eq:eta} and the cocycle identity imply
\begin{equation}\label{eq:eta-fixed}
\eta(\rho_\ep(\om),x)=\ga(\rho_\ep(\om),x)
\qquad(x\notin B_\ep).
\end{equation}
No explicit formula is needed on $B_\ep$.

Choose a bounded homogeneous cocycle $c$ such that
$\zeta=[c]\in EH_b^n(G_g)$ maps to
$\xi\in\nbar{EH}_b^n(G_g)$ under the reduced quotient. For
$\bar\om=(\om_0,\ldots,\om_n)$, the cochain representing the pullback
of the transfer is
\[
\bar\om\longmapsto
\int_{N_g}c\bigl(\eta(\rho_\ep(\om_0),x),\ldots,
                  \eta(\rho_\ep(\om_n),x)\bigr)\,d\mu(x).
\]
The integral over each of the invariant pieces above is a homogeneous
cocycle, by the proof of Proposition~\ref{p:I_b} applied to that piece.

Let $c_{\ep,A}\in C_b^n(F_2)^{F_2}$ denote the contribution of
$A_\ep$ to this cocycle. Equations \eqref{eq:trajectory-pieces} and
\eqref{eq:eta-fixed} give
\[
c_{\ep,A}(\om_0,\ldots,\om_n)
=\int_{A_\ep}
 c\bigl(u_xi_*(\om_0)u_x^{-1},\ldots,
        u_xi_*(\om_n)u_x^{-1}\bigr)\,d\mu(x).
\]
We prove that $c_{\ep,A}$ is cohomologous, through a bounded primitive,
to the cocycle $\mu(A_\ep)i^*c$. This also justifies the use of varying
conjugators.
For $u\in G_g$ let $R_{u^{-1}}$ denote right multiplication by $u^{-1}$
on vertices, and define the prism operator
\[
(K_u^n c)(g_0,\ldots,g_{n-1})
=\sum_{j=0}^{n-1}(-1)^j
c(g_0,\ldots,g_j,g_ju^{-1},\ldots,g_{n-1}u^{-1}).
\]
It preserves left homogeneity, has operator norm at most $n$, independently
of $u$, and satisfies
\[
\de K_u^n c+K_u^{n+1}\de c=R_{u^{-1}}^*c-c.
\]
The identity follows by cancellation in the alternating sum; it is the
usual prism homotopy (cf.\ \cite[Section 2.1]{Hatcher}). By left
invariance of $c$, pullback by conjugation by $u$ equals
$R_{u^{-1}}^*c$. Apply this identity with $u=u_x$, pull back by $i_*$,
and integrate over $A_\ep$. More explicitly, the pointwise integral
\[
b_{\ep,A}=\int_{A_\ep}i^*\bigl(K_{u_x}^n c\bigr)\,d\mu(x)
\in C_b^{n-1}(F_2)^{F_2}
\]
is well-defined: the integrand is measurable and
$\|b_{\ep,A}\|_\infty\leq n\mu(A_\ep)\|c\|_\infty$. Since $\de c=0$,
the homotopy identity gives
\[
c_{\ep,A}-\mu(A_\ep)i^*c=\de b_{\ep,A}.
\]
Consequently the two bounded-cohomology classes are equal:
\[
[c_{\ep,A}]=\mu(A_\ep)[i^*c]
          =\mu(A_\ep)i^*(\zeta)
\quad\text{in }H_b^n(F_2)=EH_b^n(F_2).
\]
Here exactness follows from $n\geq2$ and $H^n(F_2;\mb R)=0$.
Writing $\overline{[c_{\ep,A}]}$ for the image of $[c_{\ep,A}]$ in the
reduced quotient, the corresponding equality is
\begin{equation}\label{eq:A-contribution}
\overline{[c_{\ep,A}]}=\mu(A_\ep)i^n(\xi)
\quad\text{in }\nbar{EH}_b^n(F_2).
\end{equation}

The same argument on $A_\ep^a$ and $A_\ep^b$ reduces their contributions
to pullbacks through $h_a$ and $h_b$. Both factor through an infinite
cyclic group, whose positive-degree bounded cohomology is zero.
Thus both contributions vanish in cohomology. The contribution of
$D$ also vanishes, since it factors through the trivial group.
These amenable-group vanishing statements are recalled in
\cite{Frigerio}.

Let $R_\ep$ be the cocycle obtained by integrating over $B_\ep$.
It satisfies
\[
\|R_\ep\|_\infty\leq\mu(B_\ep)\|c\|_\infty.
\]
The preceding computations give
\[
\rho_\ep^*\Ga_b^n(\zeta)
=\mu(A_\ep)i^*(\zeta)+[R_\ep]
\quad\text{in }H_b^n(F_2).
\]
As above, $EH_b^n(F_2)=H_b^n(F_2)$ for $n\geq2$, so the error class
$[R_\ep]$ is exact. Denote its image in the reduced quotient by
$\overline{[R_\ep]}$. Passing to that quotient and using
\eqref{eq:A-contribution}, we obtain
\[
\rho_\ep^n\nbar{E\Ga}_b^n(\xi)
=\mu(A_\ep)i^n(\xi)+\overline{[R_\ep]}
\quad\text{in }\nbar{EH}_b^n(F_2).
\]
The transition collars have measure tending to zero, and
\[
\mu(B_\ep)\longrightarrow0,\qquad
\mu(A_\ep)\longrightarrow\mu(N(\al)\cap N(\be))=\La.
\]
Consequently,
\begin{align*}
\bigl\|\rho_\ep^n\nbar{E\Ga}_b^n(\xi)-\La i^n(\xi)\bigr\|
&\leq|\mu(A_\ep)-\La|\,\|i^n(\xi)\|
       +\mu(B_\ep)\|c\|_\infty\\
&\longrightarrow0.
\end{align*}
This proves the lemma.
\end{proof}

\subsection{Completion of the proof}
\begin{proof}[Proof of Theorem~\ref{T:main}]
Fix $g\geq3$ and $n\geq2$. By Proposition~\ref{p:isometric-section},
the map $\overline J_g^n$ is an isometric linear section of $i^n$.
We claim that
\[
\mathcal F_g^n=
\nbar{E\Ga}_b^n\circ\overline J_g^n:
\nbar{EH}_b^n(F_2)\longrightarrow
\nbar{EH}_b^n(\mathcal H_g)
\]
is injective. If $\mathcal F_g^n(u)=0$, apply Lemma~\ref{L:main} to
$\xi=\overline J_g^n(u)$. Since $i^n(\xi)=u$, it gives
\[
\|\La u\|
=\bigl\|\rho_\ep^n\mathcal F_g^n(u)-\La u\bigr\|
\longrightarrow0.
\]
The left side is independent of $\ep$, and $\La>0$, so $u=0$ in the
reduced space. This proves injectivity and hence
\[
\dim H_b^n(\mathcal H_g)\geq
\dim\nbar{EH}_b^n(\mathcal H_g)\geq
\dim\nbar{EH}_b^n(F_2).
\]
Thus, more generally, for every $n\geq2$ we have constructed an
injection
\[
\nbar{EH}_b^n(F_2)\hookrightarrow\nbar{EH}_b^n(\mathcal H_g).
\]
In degrees two and three the source has dimension $2^{\aleph_0}$, by the
references given at the beginning of this section. In particular, both
bounded cohomology groups in the theorem are infinite-dimensional.
\end{proof}

\begin{re}
Surjectivity alone is sufficient for the dimension argument. Indeed,
Lemma~\ref{L:main} implies
\[
\ker(\nbar{E\Ga}_b^n)\subseteq\ker(i^n).
\]
Thus Corollary~\ref{c:restriction-ge4} also proves the result for
$g\geq4$ without using the section $\overline J_g^n$.
For the elementary retraction in Section~\ref{sec:genus-ge5}, the directions
are $j^n$ surjective and $r^n$ injective, with $j^nr^n=\operatorname{Id}$.
The uniform section in Proposition~\ref{p:isometric-section} avoids
separate final arguments for the three genus ranges.
\end{re}

\begin{re}
An interesting problem is to compute $H_b^n(\op{Homeo}_0(N_g,\mu))$
for $n\geq2$ and $g=1,2$. New ideas are required for this method:
$\pi_1(N_1)$ is finite and $\pi_1(N_2)$ is virtually abelian, so both
are amenable and have trivial bounded cohomology in positive degrees.
\end{re}

\begin{re}\label{r:degree-four}
It was recently proved that $H_b^4(F_2)$ is nontrivial; see
\cite{Kastenholz}. Since $H^4(F_2;\mb R)=0$, one has
$H_b^4(F_2)=EH_b^4(F_2)$. What is not currently known is whether a
nonzero class survives after quotient by the zero-seminorm classes,
that is, whether $\nbar{EH}_b^4(F_2)$ is nontrivial. Any such nonzero
reduced class would, by the argument above, give a nonzero class in
$H_b^4(\op{Homeo}_0(N_g,\mu))$ for every $g\geq3$.
\end{re}

\enlargethispage{4\baselineskip}
\bibliographystyle{plain}
\bibliography{bibliography}

\begin{thebibliography}{10}

\bibitem{Banyaga}
Augustin Banyaga.
\newblock Formes-volume sur les vari{\'e}t{\'e}s {\`a} bord.
\newblock {\em Enseign. Math. (2)}, 20:127--131, 1974.

\bibitem{Bowditch}
Brian~H. Bowditch.
\newblock Relatively hyperbolic groups.
\newblock {\em Internat. J. Algebra Comput.}, 22(3):1250016, 66, 2012.

\bibitem{Entropy}
Michael Brandenbursky and Micha{\l} Marcinkowski.
\newblock Entropy and quasimorphisms.
\newblock {\em J. Mod. Dyn.}, 15:143--163, 2019.

\bibitem{Bounded}
Michael Brandenbursky and Micha{\l} Marcinkowski.
\newblock Bounded cohomology of transformation groups.
\newblock {\em Math. Ann.}, 382(3-4):1181--1197, 2022.

\bibitem{Volume}
Michael Brandenbursky and Micha{\l} Marcinkowski.
\newblock Volume and {E}uler classes in bounded cohomology of transformation groups.
\newblock {\em Glasg. Math. J.}, 67(1):34--49, 2025.

\bibitem{Bridson}
Martin~R. Bridson and Andr{\'e} Haefliger.
\newblock {\em Metric Spaces of Non-Positive Curvature}, volume 319 of {\em Grundlehren der mathematischen Wissenschaften}.
\newblock Springer-Verlag, Berlin, 1999.

\bibitem{Brooks}
Robert Brooks.
\newblock Some remarks on bounded cohomology.
\newblock In {\em Riemann Surfaces and Related Topics: Proceedings of the 1978 Stony Brook Conference}, volume~97 of {\em Ann. of Math. Stud.}, pages 53--63. Princeton Univ. Press, Princeton, NJ, 1981.

\bibitem{BBFIPP}
Michelle Bucher, Marc Burger, Roberto Frigerio, Alessandra Iozzi, Cristina Pagliantini, and Maria~Beatrice Pozzetti.
\newblock Isometric embeddings in bounded cohomology.
\newblock {\em J. Topol. Anal.}, 6(1):1--25, 2014.

\bibitem{Hyperbolically}
F.~Dahmani, V.~Guirardel, and D.~Osin.
\newblock Hyperbolically embedded subgroups and rotating families in groups acting on hyperbolic spaces.
\newblock {\em Mem. Amer. Math. Soc.}, 245(1156):v+152, 2017.

\bibitem{Harpe}
Pierre de~la Harpe.
\newblock {\em Topics in Geometric Group Theory}.
\newblock Chicago Lectures in Mathematics. University of Chicago Press, Chicago, IL, 2000.

\bibitem{Fathi}
Albert Fathi.
\newblock Structure of the group of homeomorphisms preserving a good measure on a compact manifold.
\newblock {\em Ann. Sci. {\'E}cole Norm. Sup. (4)}, 13(1):45--93, 1980.

\bibitem{Sisto}
R.~Frigerio, M.~B. Pozzetti, and A.~Sisto.
\newblock Extending higher-dimensional quasi-cocycles.
\newblock {\em J. Topol.}, 8(4):1123--1155, 2015.

\bibitem{Frigerio}
Roberto Frigerio.
\newblock {\em Bounded Cohomology of Discrete Groups}, volume 227 of {\em Mathematical Surveys and Monographs}.
\newblock American Mathematical Society, Providence, RI, 2017.

\bibitem{Gambaudo}
Jean-Marc Gambaudo and {\'E}tienne Ghys.
\newblock Commutators and diffeomorphisms of surfaces.
\newblock {\em Ergodic Theory Dynam. Systems}, 24(5):1591--1617, 2004.

\bibitem{Gromov}
Michael Gromov.
\newblock Volume and bounded cohomology.
\newblock {\em Inst. Hautes {\'E}tudes Sci. Publ. Math.}, 56:5--99, 1982.

\bibitem{Hatcher}
Allen Hatcher.
\newblock {\em Algebraic Topology}.
\newblock Cambridge University Press, Cambridge, 2002.

\bibitem{Hyperbolicity}
Ilya Kapovich and Kasra Rafi.
\newblock On hyperbolicity of free splitting and free factor complexes.
\newblock {\em Groups Geom. Dyn.}, 8(2):391--414, 2014.

\bibitem{Kastenholz}
Thorben Kastenholz.
\newblock Non-vanishing of the fourth bounded cohomology of free groups and codimension-two subspaces, 2025.
\newblock arXiv:2503.22511.

\bibitem{G-invariant}
Morimichi Kawasaki and Mitsuaki Kimura.
\newblock $\widehat{G}$-invariant quasimorphisms and symplectic geometry of surfaces.
\newblock {\em Israel J. Math.}, 247(2):845--871, 2022.

\bibitem{Kimura}
Mitsuaki Kimura.
\newblock Gambaudo--ghys construction on bounded cohomology.
\newblock {\em J. Math. Soc. Japan}, 77(1):135--152, 2025.

\bibitem{Rishi}
Rishi Kumar.
\newblock Construction of quasi-morphisms on groups of measure-preserving homeomorphisms of non-orientable surfaces.
\newblock Master's thesis, Ben Gurion University of the Negev, 2021.

\bibitem{Moser}
J{\"u}rgen Moser.
\newblock On the volume elements on a manifold.
\newblock {\em Trans. Amer. Math. Soc.}, 120:286--294, 1965.

\bibitem{Nitsche}
Martin Nitsche.
\newblock Higher-degree bounded cohomology of transformation groups, 2021.
\newblock arXiv:2105.08698.

\bibitem{Acylindrically}
D.~Osin.
\newblock Acylindrically hyperbolic groups.
\newblock {\em Trans. Amer. Math. Soc.}, 368(2):851--888, 2016.

\bibitem{Polterovich}
Leonid Polterovich.
\newblock Floer homology, dynamics and groups.
\newblock In Paul Biran, Octav Cornea, and Fran{\c{c}}ois Lalonde, editors, {\em Morse Theoretic Methods in Nonlinear Analysis and in Symplectic Topology}, volume 217 of {\em NATO Science Series II: Mathematics, Physics and Chemistry}, pages 417--438. Springer, Dordrecht, 2006.

\bibitem{Serre}
Jean-Pierre Serre.
\newblock {\em Trees}.
\newblock Springer-Verlag, Berlin, 1980.
\newblock Translated from the French by John Stillwell.

\bibitem{Soma}
Teruhiko Soma.
\newblock Bounded cohomology and topologically tame {K}leinian groups.
\newblock {\em Duke Math. J.}, 88(2):357--370, 1997.

\end{thebibliography}
\end{document}